
%
%
%
%
%
%
%
%
%
%
%
%
\documentstyle{amsppt}
\loadbold
\def\cstar{$C^*$-algebra}
\magnification=\magstep 1

\topmatter
\title Non-commutative spheres and\\
numerical quantum mechanics
\endtitle

\author William Arveson
\endauthor

\affil Department of Mathematics\\
University of California\\Berkeley CA 94720, USA
\endaffil

\date 7 July, 1991
\enddate
\thanks
Contributed to the proceedings of a NATO conference
on operator algebras, mathematical physics and low
dimensional topology, held in Istanbul, 1--5 July 1991.
This research was supported in part by
NSF grant DMS89-12362.
\endthanks
\abstract We discuss some basic issues that arise when one
attempts to model quantum mechanical systems on a computer,
and we describe the mathematical structure of the resulting
discretized cannonical commutation relations. The \cstar s
associated with the discretized $CCRs$ are the non-commutative
spheres of Bratteli, Elliott, Evans and Kishimoto.
\endabstract
\endtopmatter
\vfill\eject
\pagebreak
\document

\subheading{1. Introduction}

We discuss some basic issues associated with the problem of modelling
quantum mechanical phenomena on a computer.  Before setting out
to write code for any program of this kind, one first has to
recast the quantum mechanical Hamiltonian into a form appropriate
for doing numerical calculations; at the same time, one has to
retain the essential features of quantum mechanics (i.e., the
uncertainty principle).  Reflection on this problem shows that
one must confront a more basic issue, namely that of ``discretizing"
the canonical operators $P$, $Q$ which satisfy the commutation
relations
$$
PQ-QP = {\dsize\frac 1i}\bold 1.
$$

For simplicity, we consider one-dimensional quantum systems.  We
will show that there is a natural (perhaps one could say an
{\it inevitable}) way to discretize the canonical operators
$(P,Q)$ so that (a) the
discretized pair conforms to the basic principles of
numerical analysis, and (b) the uncertainty principle is preserved.
We will see that there is a one-parameter family $(P_\alpha, Q_\alpha)$,
$\alpha >0$, of discretized canonical operators.
The parameter $\alpha$ is the square of the numerical step size.
For each $\alpha$, $P_\alpha$ and $Q_\alpha$ are bounded self-adjoint
operators which obey an altered form
of the canonical commutation relations.

If $\alpha$ is not a rational multiple of $2\pi$ then the unital
\cstar\ generated by $\{P_\alpha, Q_\alpha\}$ turns out to be one of
the non-commutative spheres of Bratteli {et al} [1], [2], and hence is
a simple \cstar\ which has a unique tracial state. This \cstar\ contains the
discretized Hamiltonian, and that opens the possibility of computing its
spectrum.  While much remains to be learned, there are examples in which the
spectrum of the discretized Hamiltonian is totally disconnected.  In
particular,
these \cstar s  contain many projections.  We also show how each of these
\cstar s arises as the enveloping \cstar\ associated with an appropriately
discretized version of the Weyl canonical commutation relations.  Thus, it
seems appropriate to conclude that {\it non-commutative spheres  will arise
naturally in any  serious attempt to model quantum phenomena via computer}.

It is unlikely that there is a direct connection between these
non-commutative spheres and the quantum spheres of Woronowicz and
Podle\`s [6].  I should also point out that this work is in progress,
and that full details will appear elsewhere.  I would like to thank
Paul Chernoff and Alan Weinstein for objecting strenuously to an
earlier version of this paper; their constructive criticism has
caused me to think more carefully about the exposition, with the
result that the paper is now better than it was.  Whatever
deficiencies that remain are the sole responsibility of the author.

\subheading{2. Background}

In order to tie the discussion to a specific issue, we digress briefly
in order to describe the kind of computing problem which has motivated
this work.

Consider the one-dimensional anharmonic oscillator.  Apart from a choice
of units, the relevant equation of motion is
$$
\ddot x + x + gx^3 = 0.\leqno{2.1}
$$
Here, $g$ is a nonnegative constant and $x = x(t)$ is a real-valued function
of time $t$.  In physical terms, one has a unit mass attached to one end
of a nonlinear spring, the other end being held fixed at the origin.  When
the spring is stretched to a length $L$, the restoring force is
$L + gL^3$.  The spring is linear (i.e., it obeys Hooke's law) when
$g=0$.

In order to quantize this system, one passes to the Hilbert space
$L^2(\Bbb R)$ in the Schr\"odinger picture.  Given an initial wave
function $f \in L^2(\Bbb R)$ at time zero, basic principles dictate that
$f$ evolves in time according to
$$
f(t) = e^{itH}f, \qquad\qquad t\in \Bbb R,
$$
where $H$ is the quantum mechanical Hamiltonian.  In the case of the
anharmonic oscillator, $H$ is the operator
$$
{\dsize\frac 12} P^2 + v(Q)
$$
where $P = -i d/dx$, $Q$ is multiplication by $x$, and
$v(x) = x^2/2 + gx^4/4$ is the potential function associated with the
differential equation (2.1).

It is not hard to ``draw" solutions of the classical equation (2.1)
using a modern desktop computer.  I have done this with good results for
systems of up to 50 interacting anharmonic oscillators.  For example,
even with small systems (four masses, say) one can observe interesting
chaotic behavior resulting from the nonlinearity of the system for
$g > 0$.  On the other hand, it is less clear how one might use a
computer to create an interesting graphical representation of the
quantized system. What I have in mind is to draw pictures of
a quantized system of interacting anharmonic oscillators
which resemble the pictures described above, at least
insofar as that is possible.

With every quantum system there is an associated stationary stochastic
process whose state space is the configuration space of the classical
system ([4], [5]).  Without going into the details of how this process
is constructed (see section 5), let me just say that my intention is to
construct ``typical" sample paths for this quantum process and draw them
on a computer screen.  While I am not yet finished writing the code for
this program, I am nevertheless convinced by the results below that it is
{\it possible} to do it.  Practical obstructions, relating
to the runtime for such programs, remain to be seen.

\subheading{3. Discretizing $P$ and $Q$}

We turn now to the specific issues that arise when one seeks to build
numerical approximations to a quantum system.  One must first replace
the Hamiltonian
$$
H = {\dsize\frac 12} P^2 + v(Q)
$$
with an approximation to it which is appropriate for doing numerical
studies.  Let us first consider the momentum operator $P = -id/dx$.
Letting $V$ be the one-parameter unitary group of translations
$$
V_tf(x) = f(x-t),
$$
then of course $P$ is the infinitesimal generator of $V$
$$
V_t = e^{itP}.
$$
Equivalently, we have
$$
P = \lim_{t\to 0}{\dsize\frac 1{it}}(V_t - \bold 1)
$$
in an appropriate sense on the domain of $P$.  Conventional wisdom
indicates that $P$ should be replaced with a difference operator, the
simplest one being
$$
P_\tau={\dsize\frac 1{i\tau}}(V_\tau - \bold 1).
$$
Here, $\tau$ is the numerical step size.
While this operator is certainly bounded, it is the wrong choice
for two reasons.
It is not self adjoint, and for small $\tau$ it approximates $P$ only
to terms of order $O(\tau)$, which is not very good.  Thus it is better
to use the real part of the above operator, namely
$$
P_\tau={\dsize\frac 1{2i\tau}}(V_\tau - V_{-\tau}).
$$
This is a bounded self adjoint operator which approximates $P$ (on
smooth test functions) to terms of order $O(\tau^2)$.  Using the fact
that $V_\tau = e^{i\tau P}$ we have the equivalent formula
$$
P_\tau={\dsize\frac 1\tau}\sin (\tau P).\leqno{3.1}
$$
Considering that the function
$$
f(x) = {\dsize\frac {\sin (\tau x)}\tau}
$$
approximates the function $g(x)=x$ very closely in a fixed interval
$-M \leq x\leq +M$ when $\tau$ is small, one can expect good numerical
results from $P_\tau$ for small values of $\tau$.

We must now approximate $Q$ in a similar fashion, but in doing so we
must be careful to preserve the essential features of quantum mechanics.
That means the uncertainty principle.  In fact, the uncertainty
principle for the approximating pair $(P_\tau,Q_\tau)$ is achieved
by choosing
$$
Q_\tau = F^{-1}P_\tau F,
$$
where F is the Fourier transform operator on $L^2(\Bbb R)$.  In order
to discuss our reasons for making this choice, we digress momentarily
to a more coordinate-free setting.

Suppose that one has a pair of one-parameter unitary groups $U$, $V$
which act on the same Hilbert space $H$ in such a way that
$V_sU_t = e^{ist}U_tV_s$, and such that the set of operators
$\{U_t\}\cup\{V_s\}$ is irreducible.  Then the self-adjoint generators
$P$, $Q$ of these two groups
$$
V_t = e^{itP},\qquad U_t = e^{itQ}
$$
obey the commutation relation cited in the introduction.  Choose
a unit vector $f$ in the underlying Hilbert space $H$.  Then the probability
measure associated with the observable $Q$ and the state $f$ is the
probability measure defined (on the Borel subsets of $\Bbb R$) by
$$
q(S) = <E_Q(S)f,f>, \qquad S\subseteq\Bbb R
$$
where $E_Q$ is the spectral measure of $Q$.  Similarly, the probability
distribution associated with $P$ and $f$ is
$$
p(S) = <E_P(S)f,f>, \qquad S\subseteq\Bbb R.
$$
The uncertainty principle for the pair of observables $(P,Q)$
makes the following qualitative assertion: if a wave function $f$
represents a physical state such that the measure $q$ appears
highly concentrated near a point, then the measure $p$ must
appear ``spread out".  A similar assertion applies, of course, with
$p$ and $q$ interchanged.

It is not hard to reformulate this somewhat vague assertion entirely
in terms of the Fourier transform operator.  To see this, notice first
that there is a unitary operator $F\in \Cal B(H)$ such that
$$
P = F^{-1}QF,\qquad -Q = F^{-1}PF,\leqno{3.2}
$$
Indeed, since the pair $(P',Q')$ defined by
$P' = -Q$, $Q' = P$ is irreducible and
satisfies the same commutation relation
(equivalently, $U_{-s}V_t = e^{ist}V_tU_{-s}$, $s,t\in \Bbb R$),
the Stone-von Neumann theorem implies that the two pairs $(P,Q)$
and $(P',Q')$ are unitarily equivalent, and $F$ may be chosen as
a unitary operator implementing this equivalence.  Moreover,
because of the irreducibility hypothesis $F$ is uniquely defined
by (3.2) up to a scalar multiple of absolute value $1$.
It follows from the relation $P = F^{-1}QF$ and elementary
functional calculus that
$$
E_P(S) = F^{-1}E_Q(S)F, \qquad S\subseteq\Bbb R,
$$
and hence the formula for the measures $p$ and $q$ can be rewritten
in terms of $F$ as follows
$$
\align
q(S) =& <E_Q(S)f,f>\\
p(S) =& <E_Q(S)Ff,Ff>, \qquad S\subseteq\Bbb R.
\endalign
$$

Applying the Stone-von Neumann theorem once again, we may assume
(after replacing $(P,Q)$ with a unitarily equivalent pair)
that $H = L^2(\Bbb R)$, $P = -id/dx$, $Q = \,\,${\it multiplication by }$ x$,
and that $f$ is a normalized function in $L^2(\Bbb R)$.
Notice too that, after this change of coordinates, $F$ must be a
scalar times the usual Fourier transform operator.  Indeed, for
this choice of $P$ and $Q$ acting on $L^2(\Bbb R)$, the
Fourier transform operator $F_0$ is well-known to satisfy (3.2), i.e.,
$$
F_0P = QF_0, \qquad F_0Q = -PF_0.
$$
Hence by irreducibility we must have $F = \lambda F_0$ where
$\lambda$ is a complex number of absolute value $1$.  Thus in
these coordinates the preceding formulas for $p$ and $q$ become
$$
\align
q(S) =& \int_{S} |f(x)|^2\, dx\\
p(S) =& \int_{S} |\hat f(x)|^2\, dx, \qquad S\subseteq \Bbb R
\endalign
$$
$\hat f$ denoting the Fourier transform of the function $f$.

At this point, one can recognize the uncertainty principle as
reflecting a familiar property of the Fourier transform.
Indeed, if the measure $q$ is
highly concentrated near a point $x_0$ then so is the wave function
$f$; hence its Fourier transform $\hat f$ is highly spread out
(approximating a constant multiple of the pure exponential
$u(\lambda) = \exp(ix_0\lambda$) and thus the measure $p$ is
highly dispersed.

More generally, if one starts with {\it any} bounded self-adjoint
operator $P_1$ in the von Neumann algebra generated by all
translation operators on $L^2(\Bbb R)$, then one can expect
some form of the uncertainty principle for the pair of operators
$(P_1,FP_1F^{-1})$, at least  for bounded operators $P_1$ which
are reasonable approximations to $P$.  In more detail, we first find a
bounded real valued Borel measurable function $\phi$ such that
$P_1 = \phi(P)$.  Standard functional calculus implies that
$Q_1 = FP_1F^{-1} = F\phi(P)F^{-1} = \phi(Q)$.  Thus if $f$ is
any normalized function in $L^2(\Bbb R)$ then we have
$$
\align
<E_{Q_1}(S)f,f) =& <E_Q(\phi^{-1}S)f,f>, and\\
<E_{P_1}(S)f,f) =& <E_P(\phi^{-1}S)f,f> = <E_Q(\phi^{-1}S)Ff,Ff>
\endalign
$$
for any Borel set in $\Bbb R$, $\phi^{-1}S$ denoting the inverse
image of the set $S$ under $\phi$.  It follows that the probability
measures associated with the observables $Q_1$, $P_1$ and the state
defined by $f$ are given by
$$
\align
q_1(S) =& \int_{\phi^{-1}S} |f(x)|^2\, dx\\
p_1(S) =& \int_{\phi^{-1}S} |\hat f(x)|^2, dx, \qquad S\subseteq \Bbb R
\endalign
$$

In the most degenerate cases where $\phi$ is nearly constant, these
formulas have little content.  However, if $|\phi(x) - x|$ is small
for all $x$ in a very large interval, then these two measures closely
approximate the measures $p$, $q$ of the preceding paragraphs, and
hence they can be said to obey the uncertainty principle.

Returning now to the problem of approximating $Q$, we
conclude from these observations that in order to preserve the
uncertainty principle for the approximating operators $(P_\tau,Q_\tau)$,
one should set
$$
Q_\tau = F P_\tau F^{-1},
$$
where $F$ is the usual Fourier transform.
Notice that {\it there is no arbitrariness in the choice of $Q_\tau$};
once we have settled on a discretized approximant $P_\tau$ to $P$, $Q_\tau$
is forced upon us by requiring that the uncertainty principle should
be preserved.

Equivalently, the definition of $Q_\tau$ is
$$
Q_\tau = \frac{\sin (\tau Q)}\tau.\leqno{3.3}
$$
In terms of the one-parameter group of multiplication operators
$$
U_tf(x) = e^{itx}f(x),
$$
this definition can also be written
$$
Q_\tau = {\dsize\frac 1{2i\tau}}(U_\tau - U_{-\tau}),
$$
$\tau$ being the same step size used for $P_\tau$.

Thus for every $\tau > 0$, we have a pair of bounded self adjoint operators
$(P_\tau,Q_\tau)$.  We will write $C^*(P_\tau,Q_\tau)$ for the unital \cstar\
generated by
$P_\tau$ and $Q_\tau$.

Finally, our approximation to the original Hamiltonian is defined as
$$
H = {\dsize\frac 12} P_\tau^2 + v(Q_\tau).
$$
We abuse notation slightly here, but since we will have no need to
refer back to the original unbounded Hamiltonian it will cause no
problem. $H$ is a bounded self adjoint operator for every continuous
real-valued function $v$, and it obviously belongs to $C^*(P_\tau,Q_\tau)$.
$$
\gamma_t(a) = e^{itH}ae^{-itH}
$$
defines a uniformly continuous one-parameter group of inner automorphisms
of this \cstar.

Let us now recall the definition of non-commutative spheres $\Cal B_\alpha$
[1], [2].
Let $\alpha \in \Bbb R$ and let $U$, $V$ be a pair of unitary operators
which satisfy
$$
VU = e^{i\alpha}UV.
$$
The \cstar\ $\Cal A_\alpha$ generated by $U$ and $V$ is well known to be
independent
of the particular choice of $U$ and $V$, is simple, and has a unique
tracial state.

Noting that the pair $(U^{-1},V^{-1})$ satisfy the same
commutation relations as $(U, V)$, we can define a unique
$*$-automorphism $\sigma$ of $\Cal A_\alpha$ by requiring

$$
\sigma(U) = U^{-1}, \qquad\sigma(V) = V^{-1}.\leqno{3.4}
$$

We have $\sigma^2 = id$, and hence $\sigma$ defines a $\Bbb Z_2$-action
on $\Cal A_\alpha$.  When $\alpha$ is {\it not} a rational multiple
of $\pi$, $\Cal B_\alpha$ is defined as the fixed subalgebra

$$
\Cal B_\alpha = \{a\in \Cal A_\alpha : \sigma(a) = a \}.
$$

\proclaim{Theorem 3.5}For every positive $\tau$ such that
$\tau^2/\pi$ is irrational, $C^*(P_\tau,Q_\tau)$ is isomorphic to
$\Cal B_{\tau^2}$.
\endproclaim

\remark{Remarks}In case $\alpha/\pi$ is irrational, it is well known
that the unique trace on $\Cal A_\alpha$  gives rise to a representation of
$\Cal A_\alpha$ which generates the hyperfinite $II_1$ factor $R$, and one
can show that the weak closure of the $\sigma$-fixed subalgebra is a
subfactor of $R$ having Jones index 2.  Recalling  that such subfactors of
$R$ are unique up to conjugacy (see [3], pp.1--2 for more detail), we can
recognize this as a very stable  invariant for the embedding of
$C^*(P_\tau,Q_\tau)$ in $\Cal A_{\tau^2}$ at the  level of von Neumann
algebras.  Thus, we may think of  $C^*(P_\tau,Q_\tau)$ as being  isomorphic to
an
``index two" $C^*$-subalgebra  of the irrational  rotation algebra $\Cal
A_{\tau^2}$.

We also remark that since the step size $\tau$ must be chosen to be a small
rational number in any realistic computational setting (for example,
$\tau = 10^{-4}$ is typical for me), the \cstar s $C^*(P_\tau,Q_\tau)$ that
actually
arise will always satisfy the hypothesis of Theorem 3.5.

Let $tr$ be the tracial state of $C^*(P_\tau,Q_\tau)$.  Then for every
$\beta > 0$ there is a KMS state $\omega_\beta$ associated with the
group $\gamma$ and the value $\beta$.  $\omega_\beta$ is defined
on $C^*(P_\tau,Q_\tau)$ by
$$
\omega_\beta(a) = \frac{tr(e^{-\beta H}a)}{tr(e^{-\beta H})}.
$$
The GNS construction applied to $\omega_\beta$ leads to a representation
of $C^*(P_\tau,Q_\tau)$ which generates the hyperfinite $II_1$ factor $R$.
We shall return to the state $\omega_\beta$ in section 5.
\endremark

\subheading{4. The discretized $CCR$ algebra.}

We now give an alternate description of $C^*(P_\tau,Q_\tau)$ as the universal
\cstar\ associated with a  discretized version of Weyl's form of the
canonical commutation relations.   Significantly,
it is this realization of $C^*(P_\tau,Q_\tau)$ that is most useful for doing
numerical calculations. Fix $\alpha > 0$, not a rational multiple of
$\pi$.  Consider the discrete abelian group $G = \Bbb Z\times\Bbb Z$,
and let $\omega : G\times G \to \Bbb T$ be the bi-character defined by
$$
\omega((m,n), (p,q)) = e^{i(np - mq)\alpha/2}.
$$
By a representation of the {\it discretized canonical commutation
relations} we mean a uniformly bounded family $\{D_x: x \in G\}$ of
bounded {\it self adjoint} operators on a Hilbert space which satisfy
$$
D_xD_y = \omega(x,y)D_{x+y} + \omega(y,x)D_{x-y}.\leqno{4.1}
$$

\remark{Remarks}Formula 4.1 is a generalization of the elementary
trigonometric identity
$$
2\cos A\cos B = \cos(A+B) + \cos(A-B),
$$
in which phase shifts have been added via the cocycle $\omega$.
For example, the real-valued function $D : G\to\Bbb R$
defined (for $\sigma$ and $\phi$ fixed) by
$$
D(m,n) = 2\cos(m\sigma + n\phi)
$$
satisfies 4.1 for the trivial cocycle $\omega = 1$.
It is also clear that 4.1 is closely related to Weyl's
formulation of the canonical commutation relations.
Finally, it can be shown that any bounded family of
self-adjoint operators satisfying 4.1 must in fact obey
$\|D_x\| \leq 2$; thus {\it there is a natural
construction which produces a universal enveloping \cstar\
$C^*(G,\alpha)$ for the relations 4.1}.

It can also be shown that the original operators $P_\tau$, $Q_\tau$ serve to
define a unique $*$-representation $\pi$ of $C^*(G,\tau^2)$.  $\pi$ is
specified uniquely by requiring that
$$
\pi(D(1,0)) = 2\tau P_\tau, \qquad \pi(D(0,1)) = 2\tau Q_\tau.\leqno{4.2}
$$
While a fairly explicit formula can be given for each of the
operators $\pi(D(m,n))$ in terms of $P_\tau$ and $Q_\tau$, we shall not
do so here.  Now we are in position to formulate the second
characterization of $C^*(P_\tau,Q_\tau)$.
\endremark

\proclaim{Theorem 4.3}The representation
$\pi$ defined by 4.2 gives rise to an isomorphism of \cstar s
$$
C^*(\Bbb Z\times\Bbb Z, \tau^2) \cong C^*(P_\tau,Q_\tau).
$$
\endproclaim

\subheading{5. The canonical stochastic process}

We now indicate how one can construct a stochastic process using the
discretized operators $P_\tau$, $Q_\tau$ from section 3.  Let $v$ be
any continuous real-valued potential function, and set
$$
H = {\dsize\frac 12} P_\tau^2 + v(Q_\tau).
$$

\proclaim{Theorem 5.1}Let $f_0, f_1, \dots,f_n$ be nonnegative
continuous functions of a real variable, and let
$t_1, t_2, \dots, t_n \geq 0$.  Let $tr$ be the tracial state on
$C^*(P_\tau,Q_\tau)$.  Then
$$
tr(f_0(Q_\tau)e^{-t_1H}f_1(Q_\tau)e^{-t_2H}\dots e^{-t_nH}f_n(Q_\tau))\geq 0.
$$
\endproclaim

Using Theorem 5.1, we can make use of apparatus developed by
Klein and Landau [4] as follows.  Let $\alpha \geq 0$ and let
$tr$ be the tracial state of Theorem 5.1.

\proclaim{Theorem 5.2}Fix $\beta > 0$.  There is a real-valued
stochastic process $\{X_t: 0 \leq t \leq \beta\}$ which is defined
by the conditions
$$
\displaylines{\quad
E(f_1(X_{t_1})f_2(X_{t_2})\dots f_n(X_{t_n})) =\hfill\cr\hfill
\frac
{tr(e^{-t_1H}f_1(Q_\tau)e^{-(t_2-t_1)H}f_2(Q_\tau)\dots
																				e^{-(t_n-t_{n-1})H}f_n(Q_\tau))}
{tr(e^{-\beta H})}\quad\cr}
$$
for every $0\leq t_1\leq t_2\leq\dots\leq t_n = \beta$, every
$f_1,f_2,\dots,f_n\in C_0(\Bbb R)$, and every $n = 1,2,\dots$.
\endproclaim

\remark{Remarks}This process $\{X_t : 0\leq t\leq \beta \}$ is
stochastically continuous and periodic with period $\beta$ in that
$X_\beta(\omega) = X_0(\omega)$ almost surely.  Moreover, if one
extends it in the obvious way to a periodic process (with period $\beta$)
defined for all $t\in \Bbb R$, the resulting process is stationary
and obeys a form of reflection-positivity (i.e., Osterwalder-Schrader
positivity) which is appropriate for periodic processes.
\endremark

\enddocument